\documentclass[amstex]{amsart}
\usepackage{amssymb}
\usepackage{amsmath, amsthm}
\usepackage{mathrsfs}
\usepackage{amscd}
\usepackage{mathtools}
\mathtoolsset{showonlyrefs=true}

\input cyracc.def

\newtheorem{theorem}{Theorem}[section]
\newtheorem{lemma}[theorem]{Lemma}

\newtheorem{proposition}[theorem]{Proposition}
\newtheorem{corollary}[theorem]{Corollary}
\newtheorem{claim}{Claim}

\theoremstyle{definition}
\newtheorem{definition}[theorem]{Definition}

\theoremstyle{remark}
\newtheorem{remark}[theorem]{Remark}

\numberwithin{equation}{section}

\begin{document}

\title[The Teichm\"uller space of Zygmund diffeomorphisms]
{The real analytic structure of the Teichm\"uller space of 
circle diffeomorphisms with Zygmund continuous derivatives}

\author[K. Matsuzaki]{Katsuhiko Matsuzaki}
\address{Department of Mathematics, School of Education, Waseda University, 
Shinjuku, Tokyo 169-8050, Japan}
\email{matsuzak@waseda.jp}
\thanks{Research supported by 
Japan Society for the Promotion of Science (KAKENHI 23K25775 and 23K17656)}

%\dedicatory{Dedicated to Professor Yunping Jiang on the occasion of his 65th birthday}
 

\begin{abstract}
We apply the methods of simultaneous uniformization and composition operators on
Besov spaces to the Teichm\"uller space $T^Z$ of 
circle diffeomorphisms with Zygmund continuous deriv\-atives.
As consequences, we obtain the following: (1) a new proof of the correspondence between quasiconformal self-homeomorphisms of
the unit disk with complex dilatations of linear decay order and their quasisymmetric extensions to the unit circle
with regularity 
in the Zygmund continuously differentiable class; (2) a real-analytic equivalence of $T^Z$ with
the real Banach space of Zygmund continuous functions on the unit circle.
\end{abstract}

\maketitle

%\vspace{-0.7cm}

\section{Introduction}

Let ${\rm Diff}^{1+Z}(\mathbb S)$ be the group of all orientation-preserving diffeomorphisms 
$h:\mathbb S \to \mathbb S$ of the unit circle whose non-degenerate derivatives $h'$ are continuous and satisfy the {\it Zygmund condition}
\begin{equation}\label{zy}
|h'(e^{i(\theta+t)}) - 2h'(e^{i\theta}) + h'(e^{i(\theta-t)})| \leq Ct
\end{equation}
for all $e^{i\theta} \in \mathbb S$ and $t > 0$. Here, $C>0$ is a constant independent of $e^{i\theta}$ and $t$.
A continuous function $h'$ satisfying this condition has the modulus of continuity
\begin{equation}\label{modulus}
|h'(e^{i(\theta+t)}) - h'(e^{i\theta})| = O(t\log(1/t)) \quad (t \to 0).
\end{equation}
It follows that $h'$ is $\alpha$-H\"older continuous for every $\alpha \in (0,1)$; namely,
$|h'(e^{i(\theta+t)}) -h'(e^{i\theta})| \leq C't^\alpha$ holds.
The group of all orientation-preserv\-ing diffeomorphisms 
$h$ of $\mathbb S$ such that $h'$ is $\alpha$-H\"older continuous is denoted by ${\rm Diff}^{1+\alpha}(\mathbb S)$.
For these elements, the normalization is imposed by fixing three points $1$, $i$, and $-i$,
and the subgroups consisting of all normalized elements are denoted by 
${\rm Diff}^{1+Z}_*(\mathbb S)$ and ${\rm Diff}^{1+\alpha}_*(\mathbb S)$, respectively.
As we will see below, ${\rm Diff}^{1+Z}_*(\mathbb S)$ can be regarded as 
the Teichm\"uller space $T^Z$, and ${\rm Diff}^{1+\alpha}_*(\mathbb S)$ as $T^\alpha$.
In this paper, however, we focus on $T^Z$, which is defined as a subspace of
the universal Teichm\"uller space $T$ in the following way.

The space of Beltrami coefficients on the exterior unit disk
$\mathbb D^* = \{z \mid {|z|>1}\} \cup \{\infty\}$ is defined by
$$
M(\mathbb D^*) = \{\mu \in L_\infty(\mathbb D^*) \mid \Vert \mu \Vert_\infty < 1\}.
$$
Let $H(\mu):\mathbb D^* \to \mathbb D^*$ be the normalized quasiconformal self-homeo\-mor\-phism 
whose complex dilatation $\bar \partial H/\partial H$ is $\mu \in M(\mathbb D^*)$.
Such an $H(\mu)$ extends uniquely to the unit circle $\mathbb S$ 
as a quasisymmetric self-homeomorphism $h_\mu$.
The normalization of $H(\mu)$ is given by that of $h_\mu$, and hence $H(\mu)$ is uniquely determined by $\mu \in M(\mathbb D^*)$.
Denote by ${\rm QS}_*(\mathbb S)$ the group of all normalized quasisymmetric self-homeomorphisms of $\mathbb S$.
The following inclusion relations hold: 
$$
{\rm Diff}^{1+Z}_*(\mathbb S) \subset {\rm Diff}^{1+\alpha}_*(\mathbb S) \subset {\rm QS}_*(\mathbb S).
$$

The boundary extension $H(\mu)|_{\mathbb S}=h_\mu$ 
defines a surjection
$\pi:M(\mathbb D^*) \to {\rm QS}_*(\mathbb S)$. 
The Teichm\"uller equivalence $\mu \sim \nu$ in $M(\mathbb D^*)$ is defined by the coincidence
$h_\mu = h_\nu$ in ${\rm QS}_*(\mathbb S)$, and
the universal Teich\-m\"uller space $T$ is defined to be
the quotient space $M(\mathbb D^*)/\sim$. Hence, $T$ can be identified with ${\rm QS}_*(\mathbb S)$.
We call this quotient map and also $\pi$ the {\it Teichm\"uller projection}.

To introduce the Teichm\"uller space $T^Z$,
we impose the linear degeneration condition at the boundary on Beltrami coefficients:
\begin{align*}
M^Z(\mathbb D^*) = \{ \mu \in M(\mathbb D^*) \mid \Vert \mu \Vert_Z
= \underset{|z|>1}{\mathrm{ess\,sup}}\, ((|z|-1)^{-1} \lor 1) |\mu(z)| < \infty \}.
\end{align*}
Here, the symbol $\lor$ stands for the maximum of the two values. 
Tang and Wu \cite[Theorem 1.1]{TW1} proved 
the precise correspondence of $M^Z(\mathbb D^*)$ to ${\rm Diff}^{1+Z}_*(\mathbb S)$
under the Teichm\"uller projection $\pi:M(\mathbb D^*) \to {\rm QS}_*(\mathbb S)$.

\begin{theorem}\label{basic}
$\pi(M^Z(\mathbb D^*)) = {\rm Diff}^{1+Z}_*(\mathbb S)$. 
\end{theorem}

Thus, the Teichm\"uller space $T^Z$, defined as $M^Z(\mathbb D^*)/\sim$, is 
identified with ${\rm Diff}^{1+Z}_*(\mathbb S)$.
In the first part of this paper, we reconstruct the proof of this theorem.
We note that the Teichm\"uller space $T^\alpha$ $(0<\alpha <1)$, defined as $M^\alpha(\mathbb D^*)/\sim$, is 
identified with ${\rm Diff}^{1+\alpha}_*(\mathbb S)$ in \cite{M0}, where 
$M^\alpha(\mathbb D^*)$ consists of all Beltrami coefficients $\mu \in M(\mathbb D^*)$
satisfying $\Vert \mu \Vert_{(\alpha)}<\infty$. See Theorem \ref{alpha2} for the definition of this norm.

To clarify the arguments, we divide the statement of Theorem \ref{basic} into two parts:

\begin{claim}\label{claim1}
$\pi(M^Z(\mathbb D^*))$ is contained in ${\rm Diff}^{1+Z}_*(\mathbb S)$.
\end{claim}

\begin{claim}\label{claim2}
$\pi:M^Z(\mathbb D^*) \to {\rm Diff}^{1+Z}_*(\mathbb S)$ is surjective; that is, for any $h \in {\rm Diff}^{1+Z}_*(\mathbb S)$, there exists a quasiconformal extension $H:\mathbb D^* \to \mathbb D^*$ whose
complex dilatation $\mu_H=\bar \partial H/\partial H$ belongs to $M^Z(\mathbb D^*)$.
\end{claim}

The theory of the Teichm\"uller space $T^Z = M^Z(\mathbb D^*)/\sim$ has been developed
from this definition involving quasiconformal mappings, similarly to the universal Teichm\"uller space (see \cite[Chapter 3]{Le}).
For $\mu \in M(\mathbb D^*)$, let $F_\mu:\mathbb D \to \mathbb C$ be the conformal map of the unit disk
$\mathbb D$ that
admits a quasiconformal extension to the Riemann sphere $\widehat{\mathbb C}$ with complex dilatation on 
$\mathbb D^*$,
satisfying the normalization $F_\mu(0) = 0$, $F_\mu'(0) = 1$, and $F_\mu(\infty) = \infty$.
To consider its pre-Schwarzian derivative $N_{F_\mu}=(\log(F_\mu)')'$ and 
the Schwarzian derivative $S_{F_\mu}=(N_{F_\mu})'-(N_{F_\mu})^2/2$, 
we prepare the corresponding Banach spaces in ${\rm Hol}(\mathbb D)$, the space of holomorphic functions on $\mathbb D$:
\begin{align*}
B^Z(\mathbb D)& = \{\Phi \in {\rm Hol}(\mathbb D) \mid \Vert \Phi \Vert_{B^Z} < \infty\},\\
&\Vert \Phi \Vert_{B^Z} = |\Phi'(0)| + \sup_{|z|<1}(1 - |z|^2)|\Phi''(z)|;\\
A^Z(\mathbb D)&= \{\Psi \in {\rm Hol}(\mathbb D) \mid \Vert \Psi \Vert_{A^Z} < \infty\},\\
&\Vert \Psi \Vert_{A^Z} = \sup_{|z|<1}(1 - |z|^2)|\Psi(z)|.
\end{align*}
Here, we regard $B^Z(\mathbb D)$ as a Banach space modulo constant functions.

The complex-analytic theory of $T^Z$ is based on the following characterization of $\mu \in M^Z(\mathbb D^*)$,
which is summarized in \cite{TW1}. After this work, a complex Banach manifold structure was
provided for $T^Z$ via the Schwarzian derivative map $S:M(\mathbb D^*) \to A^Z(\mathbb D)$ 
defined by $\mu \mapsto S_{F_\mu}$, which factors through
the Teichm\"uller projection into
the Bers embedding $\alpha:T^Z \to A^Z(\mathbb D)$ (see \cite[Theorem 3]{M3}). Moreover,
the fiber bundle $\widetilde{\mathcal T}^Z$ over $T_Z$ in $B^Z(\mathbb D)$ is induced by
the {\it pre-Schwarzian derivative map} $L:M(\mathbb D^*) \to B^Z(\mathbb D)$
defined by $\mu \mapsto \log (F_\mu)'$ (see \cite[Theorem 2]{M4}).
Both $S$ and $L$ are holomorphic split submersions onto their images.

\begin{theorem}\label{equivalence}
For a conformal homeomorphism $F_\mu:\mathbb D \to \mathbb C$ given by $\mu \in M(\mathbb D^*)$,
the following conditions are equivalent:
\begin{enumerate}
\item
$\mu \in M^Z(\mathbb D^*)$;
\item
$\log (F_\mu)' \in B^Z(\mathbb D)$;
\item
$S_{F_\mu} \in A^Z(\mathbb D)$;
\item
$(F_\mu)''' \in A^Z(\mathbb D)$;
\item
$(F_\mu)'$ extends continuously to $\mathbb S$ and $(F_\mu)'|_{\mathbb S}$ 
satisfies the Zygmund condition;
\item
$\log (F_\mu)'$ extends continuously to $\mathbb S$ and $\log (F_\mu)'|_{\mathbb S}$ 
satisfies the Zygmund condition.
\end{enumerate}
\end{theorem}

The implication $(1) \Rightarrow (2)$ follows from Dyn\cprime kin \cite[Theorem 1]{Dy}, 
whereas $(2) \Rightarrow (1)$ was proved by Becker and Pommerenke \cite[Satz 4]{BP}.
Moreover, $(3) \Rightarrow (1)$ is shown by Becker \cite[Theorem 3]{Be}.
In Theorem \ref{alpha2}, we give a new proof of $(1) \Rightarrow (3)$.
The equivalence $(4) \Leftrightarrow (5)$ is due to Zygmund \cite[Theorem 13]{Zy}; see also Theorem \ref{zygmund}. 
Tang and Wu \cite[Theorem 2.6]{TW1} proved the chain of implications $(3) \Rightarrow (4) \Rightarrow (2) \Rightarrow (3)$ in this order. 
The equivalence $(2) \Leftrightarrow (6)$ also follows from Zygmund's theorem. 
Under the condition that $(F_\mu)'|_{\mathbb S}$ does not vanish, we have $(5) \Leftrightarrow (6)$.
See the second paragraph of Section \ref{3} concerning the post-composition of the logarithm
for the Zygmund condition.

\begin{remark}
In the situation of the above theorem, $F_\mu$ maps $\mathbb D$ conformally onto a domain in $\mathbb C$
bounded by a non-degenerate $C^{1+Z}$ Jordan curve.
This condition also yields an equivalent characterization (see \cite[Theorem 2.4]{TW1}),
which is an extension of the Kellogg--Warschawski theorem (see \cite[Theorem 3.6]{Pom}).
\end{remark}

For later use, we define the
following Banach spaces (modulo constants) of Zygmund continuous and differentiable functions:
\begin{align*}
C^Z(\mathbb S)&=\{\phi \in C^0(\mathbb S) \mid \Vert \phi \Vert_{C^Z}<\infty\},\\
& \Vert \phi \Vert_{C^Z}
=\sup_{e^{i\theta} \in {\mathbb S},\, t>0}\frac{|\phi(e^{i(\theta+t)}) - 2\phi(e^{i\theta}) + \phi(e^{i(\theta-t)})|}{2t};\\
C^{1+Z}(\mathbb S)&=\{\psi \in C^1(\mathbb S) \mid \Vert \psi \Vert_{C^{1+Z}}<\infty\},\\
& \Vert \psi \Vert_{C^{1+Z}}
=\Vert \psi' \Vert_{C^{Z}}+\Vert \psi' \Vert_{L_\infty}.
\end{align*}
We can alternatively define them as the homogeneous Besov spaces $\dot B_{\infty,\infty}^1(\mathbb S)$ and 
$\dot B_{\infty,\infty}^2(\mathbb S) \cap \dot W^1_\infty(\mathbb S)$ with equivalent norms, respectively; see Section \ref{2}.

Concerning Theorem \ref{basic},
the argument of Tang and Wu in its proof relies on the pseudoanalytic extension
developed by Dyn\cprime kin \cite[Theorem 1]{Dyps}. In fact, the following statement in a special case is used; see \cite[Lemma 2.1]{TW1}.

\begin{theorem}\label{pseudoanalytic}
Let $\psi$ be a continuous function on $\mathbb S$ that extends holomorphically to $\mathbb D$.
Then, 
$\psi$ belongs to $C^{1+Z}(\mathbb S)$ 
if and only if there exists a $C^1$ extension $\Psi$ of $\psi$ to $\mathbb D^*$ such that
$$
\sup_{|z|>1} (|z|-1)^{-1}|\bar\partial \Psi(z)| < \infty.
$$
\end{theorem}

In the case of one lower order of differentiability, 
the {\it Beurling--Ahlfors extension} \cite{BA} satisfies a similar condition, 
as shown by Gardiner and Sullivan \cite[p.~733]{GS}.
This was used for describing intrinsically the tangent space of ${\rm QS}_*(\mathbb S)$ at the identity.
See also \cite{GH} and \cite{nag}.

\begin{proposition}\label{BACZ}
A function $\phi$ on $\mathbb S$ belongs to $C^Z(\mathbb S)$ 
if and only if its Beurling--Ahlfors extension $\Phi$ to $\mathbb D^*$ satisfies
$
\sup_{|z|>1} |\bar\partial \Phi(z)| < \infty.
$
\end{proposition}

However, the Beurling--Ahlfors extension does not work properly
for the purpose of Claim \ref{claim2}. 
Consider an orientation-preserving diffeo\-morphism $h:\mathbb S \to \mathbb S$ whose derivative is Lipschitz continuous:
$$
|h'(e^{i(\theta+t)}) - h'(e^{i\theta})| \leq Ct.
$$ 
This Lipschitz condition is stronger than the Zygmund condition \eqref{zy}.
The following characterization of such mappings was proved by Hu \cite[Theorem 10]{Hu},
which makes our problem rather delicate.

\begin{theorem}\label{junhu}
Let $h$ be an orientation-preserving self-diffeo\-mor\-phism of $\mathbb S$.
A necessary and sufficient condition for $h'$ to be Lipschitz continuous is that the complex dilatation $\mu_H$ 
of the Beurling--Ahlfors extension $H:\mathbb D^* \to \mathbb D^*$ of $h$ belongs to $M^Z(\mathbb D^*)$.
\end{theorem}

The Beurling--Ahlfors extension is originally
defined for quasisymmetric homeomorphisms of the real line $\mathbb R$.
For those $h$ on the unit circle $\mathbb S$, we lift $h$ to $\mathbb R$ via the universal covering $\bar p:\mathbb R \to \mathbb S$
given by $x \mapsto e^{ix}$, and apply the Beurling--Ahlfors extension to these lifts $\tilde h$ to obtain
equivariant quasiconformal self-homeomorphisms $\widetilde H$ with $\widetilde H(z+2\pi)=\widetilde H(z)+2\pi$ 
of the lower half-plane $\mathbb H^*$.
They are projected down to 
quasiconformal self-homeomorphisms of
$\mathbb D^* \setminus \{\infty\}$ by the universal covering $p:z \mapsto e^{iz}$, which extend quasiconformally to $\infty$.
See the following diagram:
\medskip
\begin{equation}\label{univ}
\begin{CD}
     \qquad \mathbb R \curvearrowleft \tilde h @>{{\rm BA \, ext.}}>> \widetilde H \curvearrowright \mathbb H^*\qquad \\
  @V{{\bar p}}VV    @VV{p}V \\
     \qquad \mathbb S \curvearrowleft h @>>> H \curvearrowright \mathbb D^* \setminus \{\infty\}.
\end{CD}
\end{equation}
\smallskip

In this paper, we provide independent solutions to Claims \ref{claim1} and \ref{claim2} by a method different from that based on pseudoanalytic extensions used in Theorem \ref{pseudoanalytic}.
Claim \ref{claim1} is addressed in Section \ref{3} using arguments of conformal welding and composition operators,
following preparations on Besov spaces and their interpolation
from real analysis in Section \ref{2}. Claim \ref{claim2} is treated in Section \ref{5} after
preparing results on boundary extensions of holomorphic functions and parametrizations of curves 
from complex analysis in Section \ref{4}. The proof of Claim \ref{claim2} is obtained by applying
the method of simultaneous uniformization. To employ this argument, we lift the relevant elements defined
on $\mathbb S$, $\mathbb D$, and $\mathbb D^*$ to periodic elements on $\mathbb R$, $\mathbb H$, and $\mathbb H^*$
via the universal cover as in diagram \eqref{univ}. 
In Section \ref{6}, as an application of the method of simultaneous uniformization, we also
present a result asserting that the Teichm\"uller space $T^Z$ is real-analytically equivalent to 
the real Banach space of Zygmund continuous functions on $\mathbb S$ (Corollary \ref{real-analytic}).
As an appendix, in Section \ref{7}, we provide a different proof of $(1) \Rightarrow (3)$ in 
Theorem \ref{equivalence}.

The methods used to prove Claims \ref{claim1} and \ref{claim2} were originally developed in \cite{M6}
for the integrable (Weil--Petersson) Teichm\"uller spaces.

\section{Besov spaces and real interpolation}\label{2}

We introduce function spaces that appear in this paper and
then show certain properties of these spaces from a functional-analytic perspective.

For $0<\alpha<1$, let $C^\alpha(\mathbb S)$ be the set of $\alpha$-H\"older 
continuous functions $\phi$ on $\mathbb S$. By lifting $\phi$ to the universal cover $\mathbb R$
and denoting it by the same symbol, we may regard $\phi$ as a periodic function on $\mathbb R$,
that is, on $\mathbb S=\mathbb R/2\pi \mathbb Z$. See diagram \eqref{univ}. Hereafter, we represent
the functions on $\mathbb S$ in this way by $x \in \mathbb R$ not using $e^{i\theta} \in \mathbb S$.
The seminorm of $C^\alpha(\mathbb S)$ is defined by
$$
\Vert \phi \Vert_{C^\alpha}=\sup_{x,y \in \mathbb S} \frac{|\phi(x)-\phi(y)|}{|x-y|^\alpha}.
$$
Let $C^L(\mathbb S)$ be the set of Lipschitz continuous functions 
$\phi$ on $\mathbb S$. 
The seminorm of $C^L(\mathbb S)$ is defined by
$$
\Vert \phi \Vert_{C^L}=\sup_{x,y \in \mathbb S} \frac{|\phi(x)-\phi(y)|}{|x-y|}.
$$
Moreover, let $C^Z(\mathbb S)$ be the set of Zygmund
continuous functions $\phi$ on $\mathbb S$. 
The seminorm of $C^Z(\mathbb S)$ is defined by
$$
\Vert \phi \Vert_{C^Z}=\sup_{x,y \in \mathbb S} \frac{|\phi(x)+\phi(y)-2\phi((x+y)/2)|}{|x-y|}.
$$

Modulo constant functions, the above seminorms become norms, and
the function spaces defined above become Banach spaces.
There are inclusion relations 
$C^L(\mathbb S) \hookrightarrow C^Z(\mathbb S) \hookrightarrow C^\alpha(\mathbb S)$,
and the inclusion maps are continuous. In particular, 
$\Vert \phi \Vert_{C^Z} \leq \Vert \phi \Vert_{C^L}$ holds.

These function spaces can be generalized by introducing Besov spaces.
For $m \in \mathbb{N}$ and $t \in \mathbb{R}$, 
the $m$-th order difference of a function $\phi$ on $\mathbb{S}=\mathbb R/2\pi \mathbb Z$ is defined inductively as
$$
\Delta_t^1 \phi(x) = \phi(x+t) - \phi(x); \quad
\Delta_t^{m+1} \phi(x) = \Delta_t^m \phi(x+t) - \Delta_t^m \phi(x).
$$

\begin{definition}
For $s >0$ and $0 < p,\, q \leq \infty$, a seminorm of an integrable 
function $\phi$ on $\mathbb{S}$ is given by
$$
\Vert \phi \Vert_{\dot B^s_{p,q}} = \left(\int_{-\pi}^{\pi}
|t|^{-sq} \Vert \Delta_t^{\lfloor s \rfloor+1} \phi \Vert_{L_p}^q \frac{dt}{|t|}
\right)^{1/q}.
$$
When $q=\infty$, the $L_q$-norm by the integration above is understood as the supremum by the usual convention.
The set of those $\phi$ with $\Vert \phi \Vert_{\dot B^s_{p,q}} < \infty$ is defined as the homogeneous {\it Besov space}
$\dot B^s_{p,q}(\mathbb S)$.
\end{definition}

In this paper, we only deal with the case $p=q=\infty$. For $0<\alpha<1$, $\dot B^\alpha_{\infty,\infty}(\mathbb S)$ is
identified with $C^\alpha(\mathbb S)$, and $\dot B^1_{\infty,\infty}(\mathbb S)$ is
identified with $C^Z(\mathbb S)$ (see \cite[p.541]{Leo}). 
By contrast, the homogeneous Sobolev space $\dot W^1_\infty(\mathbb S)$,
with seminorm $\Vert \phi \Vert_{\dot W^1_\infty}=\Vert \phi' \Vert_{L_\infty}$, is identified with $C^L(\mathbb S)$. 

We use interpolation for Besov spaces.
In general, real interpolation of Banach spaces is defined as follows.

\begin{definition}
Let $(X_0, \Vert \cdot \Vert_{X_0})$ and $(X_1, \Vert \cdot \Vert_{X_1})$ 
be Banach spaces that are continuously embedded into a common topological vector space,
which we call an admissible pair of Banach spaces. For $t>0$, 
$$
K(x,t)=\inf\,\{\Vert x_0 \Vert_{X_0}+t\Vert x_1 \Vert_{X_1} \mid x=x_0+x_1,\ x_0 \in X_0,\ x_1 \in X_1\} 
$$
is defined for $x \in X_0+X_1$. Then, for $\sigma \in (0,1)$ and $q \geq 1$, the set
$$
(X_0,X_1)_{\sigma,q}=
\{x \in X_0 +X_1 \mid \Vert x \Vert_{\sigma,q}=\left( \int_0^\infty K(x,t)^q 
\frac{dt}{t^{1+\sigma q}} \right)^{\frac{1}{q}}<\infty\}
$$ 
is defined to be the {\it real interpolation space} of $X_0$ and $X_1$.
\end{definition}

The real interpolation space $(X_0,X_1)_{\sigma,q}$ is a Banach space with the norm $\Vert \cdot \Vert_{\sigma,q}$
(see \cite[Theorem 16.5]{Leo}). 
For interpolation of Besov spaces, the following result is known (see \cite[Corollary 17.42]{Leo}).
The above definition of real interpolation can also be applied to homogeneous Besov spaces with seminorms
(see \cite[Remark 17.31]{Leo}).

\begin{proposition}\label{Besov}
Let $1 \leq p,q\leq \infty$, $0<\sigma<1$, and $0<s_0<s_1$. Then
$$
(\dot B^{s_0}_{p,q}(\mathbb S),\dot B^{s_1}_{p,q}(\mathbb S))_{\sigma,q}=\dot B^{s}_{p,q}(\mathbb S),
$$
where $s=(1-\sigma)s_0+\sigma s_1$.
\end{proposition} 

The boundedness of associated operators under real interpolation holds in the following form (see \cite[Theorem 16.12]{Leo}).

\begin{proposition}\label{interpolation}
Let $(X_0, \Vert \cdot \Vert_{X_0})$ and $(X_1, \Vert \cdot \Vert_{X_1})$ 
be an admissible pair of Banach spaces. Let $P:X_0+X_1 \to X_0+X_1$ be a linear transformation
such that its restrictions $P|_{X_0}$ and $P|_{X_1}$ are bounded linear transformations
of $X_0$ and $X_1$, respectively. Then, for any $\sigma \in (0,1)$ and $q \geq 1$,
$P|_{(X_0,X_1)_{\sigma,q}}$ is a bounded linear transformation of
the real interpolation space $(X_0,X_1)_{\sigma,q}$, and its operator norm satisfies
$$
\Vert P|_{(X_0,X_1)_{\sigma,q}} \Vert \leq \Vert P|_{X_0} \Vert^{1-\sigma} \Vert P|_{X_1} \Vert^\sigma.
$$
\end{proposition}

\section{Solution to Claim \ref{claim1}}\label{3}

Let ${\rm QS}(\mathbb S)$ denote the group of quasisymmetric 
self-homeo\-morph\-isms of $\mathbb S$. 
For $0<\alpha<1$, let ${\rm Diff}^{1+\alpha}(\mathbb S)$ be the subgroup of ${\rm QS}(\mathbb S)$
consisting of orientation-preserving self-diffeomorphisms $h:\mathbb S \to \mathbb S$ (with non-vanishing derivatives)
such that the derivative $h'$ belongs to $C^\alpha(\mathbb S)$.
This condition is equivalent to requiring $\log h' \in C^\alpha(\mathbb S)$.
The fact that ${\rm Diff}^{1+\alpha}(\mathbb S)$ is a group under composition can be verified easily.
In fact, this is a topological group in an appropriate topology (see \cite[Proposition 5.2]{M1}).

Analogously, we define ${\rm Diff}^{1+Z}(\mathbb S)$ to be the subset of ${\rm Diff}^{1+\alpha}(\mathbb S)$ consisting of all such $h$ that $h'$ belongs to $C^Z(\mathbb S)$.
In this case, this is equivalent to the condition that $\log h' \in C^Z(\mathbb S)$ since $h' \neq 0$.
Indeed, a straightforward computation of the second order difference of $\log h'$
using \eqref{modulus} yields the desired estimate.
We show later in Proposition \ref{group} that ${\rm Diff}^{1+Z}(\mathbb S)$ 
is a subgroup of ${\rm Diff}^{1+\alpha}(\mathbb S) \subset {\rm QS}(\mathbb S)$.

We utilize {\it conformal welding} to prove Claim \ref{claim1}. Namely, we decompose 
$h_\mu:\mathbb S \to \mathbb S$ for
$\mu \in M^Z(\mathbb D^*)$ into
the boundary extensions of conformal homeomorphisms $F_\mu:\mathbb D \to \Omega$ and $G^{\mu^{-1}}:\mathbb D^* \to \Omega^*$
such that $h_{\mu}=(G^{\mu^{-1}}|_{\mathbb S})^{-1} \circ F_\mu|_{\mathbb S}$. 
Here, $\mu^{-1}$ denotes the complex dilatation of $H(\mu)^{-1}$, and by assuming $\mu^{-1}$ on $\mathbb D$,
the quasiconformal homeomorphism $G^{\mu^{-1}}$ of $\widehat{\mathbb C}$ is defined by extending it conformally to $\mathbb D^*$.
From this, we obtain
\begin{equation}\label{log}
\log h_\mu' = -\log g' \circ h_\mu + \log f'
\end{equation}
for $f=F_\mu|_{\mathbb S}$ and $g=G^{\mu^{-1}}|_{\mathbb S}$. 
Note that $h_\mu$ belongs to ${\rm Diff}^{1+\alpha}(\mathbb S)$ for any $\alpha \in (0,1)$ (see \cite[Theorem 6.7]{M1}
and the references therein).

We apply Theorem \ref{equivalence} to these
conformal homeomorphisms to obtain that $\log f'$ and $\log g'$ belong to $C^Z(\mathbb S)$.
For $f$, this is straightforward. For $g$, we have to show that $\mu^{-1}$ belongs to $M^Z(\mathbb D^*)$.
For $H=H(\mu)$, we have
$$
\mu^{-1}(H(z))=-\mu(z) \frac{\partial H(z)}{\overline{\partial H(z)}} \quad (z \in \mathbb D^*).
$$
Here, by applying \cite[Theorem 6.4]{M1}, we have $|H(z)|-1 \asymp |z|-1$. Then, it follows from 
$\mu \in  M^Z(\mathbb D^*)$ that $\mu^{-1} \in  M^Z(\mathbb D^*)$.

Having $\log g' \in C^Z(\mathbb S)$, we obtain $\log g'\circ h_\mu \in C^Z(\mathbb S)$ by applying the following
claim. This is a crucial point in our proof of Claim \ref{claim1}.

\begin{lemma}\label{composition}
For any $h \in {\rm Diff}^{1+\alpha}(\mathbb S)$ with $0<\alpha<1$,
the composition operator $P_h$ on $C^Z(\mathbb S)$ defined by
$$
P_h: \phi \mapsto \phi \circ h \qquad (\phi \in C^Z(\mathbb S))
$$
is a Banach automorphism of $C^Z(\mathbb S)$. Moreover, the operator norm of $P_h$ is bounded by
$c_\alpha\Vert h \Vert_{C^{1+\alpha}}$ for some constant $c_\alpha>0$.
\end{lemma}

\begin{proof}
Recall that $C^Z(\mathbb S)$ coincides with $\dot B^1_{\infty,\infty}(\mathbb S)$.
Then, represent $\dot B^1_{\infty,\infty}(\mathbb S)$ by the real interpolation
$$
\dot B^1_{\infty,\infty}(\mathbb S)=(\dot B^{1-\alpha}_{\infty,\infty}(\mathbb S),\dot B^{1+\alpha}_{\infty,\infty}(\mathbb S))_{1/2,\infty}
$$
as a special case of Proposition \ref{Besov}.
Here, $\dot B^{1-\alpha}_{\infty,\infty}(\mathbb S)=C^{1-\alpha}(\mathbb S)$, and we set 
$$
\ddot B^{1+\alpha}_{\infty,\infty}(\mathbb S)=\dot B^{1+\alpha}_{\infty,\infty}(\mathbb S) \cap \dot W^1_\infty(\mathbb S)
$$
with seminorm $\Vert \phi \Vert_{\dot B^{1+\alpha}_{\infty,\infty}}+\Vert \phi \Vert_{\dot W^1_\infty}$. This coincides with
$$
C^{1+\alpha}(\mathbb S)=\{\psi \in C^1(\mathbb S) \mid \Vert \psi \Vert_{C^{1+\alpha}}:=
\Vert \psi' \Vert_{C^{\alpha}}+\Vert \psi' \Vert_{L_\infty}<\infty\}.
$$ 
Since $(\dot B^{1-\alpha}_{\infty,\infty}(\mathbb S),\dot W^1_\infty(\mathbb S))_{1/2,\infty}=
\dot B^{1-\alpha/2}_{\infty,\infty}(\mathbb S)$ by \cite[Theorem 17.30]{Leo},
which contains $\dot B^1_{\infty,\infty}(\mathbb S)$, we eventually have
\begin{equation}\label{represent}
C^Z(\mathbb S)=\dot B^1_{\infty,\infty}(\mathbb S)=(\dot B^{1-\alpha}_{\infty,\infty}(\mathbb S),\ddot B^{1+\alpha}_{\infty,\infty}(\mathbb S))_{1/2,\infty}.
\end{equation}

For $\phi_1 \in \dot B^{1-\alpha}_{\infty,\infty}(\mathbb S)=C^{1-\alpha}(\mathbb S)$,
we have
\begin{align}
|\phi_1 \circ h(x)-\phi_1 \circ h(y)| &\leq \Vert \phi_1 \Vert_{C^{1-\alpha}}|h(x)-h(y)|^{1-\alpha}\\
&\leq \Vert \phi_1 \Vert_{C^{1-\alpha}}\Vert h' \Vert_{L_\infty}^{1-\alpha}|x-y|^{1-\alpha}.
\end{align}
This implies that the composition operator $P_h$ is bounded on $C^{1-\alpha}(\mathbb S)$ with operator norm at most $\Vert h' \Vert_{L_\infty}^{1-\alpha} \leq \Vert h \Vert_{C^{1+\alpha}}^{1-\alpha}$.

For $\phi_2 \in \ddot B^{1+\alpha}_{\infty,\infty}(\mathbb S)=C^{1+\alpha}(\mathbb S)$,
the difference 
$$
|(\phi_2 \circ h)'(x)-(\phi_2 \circ h)'(y)|=|(\phi_2)' (h(x))h'(x)-(\phi_2)' (h(y))h'(y)|
$$
is bounded above by
\begin{align}
&\quad \ |(\phi_2)' (h(x))h'(x)-(\phi_2)' (h(y))h'(x)|\\
&\qquad \qquad \qquad \qquad \ +|(\phi_2)' (h(y))h'(x)-(\phi_2)' (h(y))h'(y)|\\
&\leq |(\phi_2)' (h(x))-(\phi_2)' (h(y))|\Vert h' \Vert_{L_\infty}+\Vert (\phi_2)' \Vert_{L_\infty}|h'(x)-h'(y)|\\
&\leq (\Vert (\phi_2)'\Vert_{C^\alpha} \Vert h' \Vert_{L_\infty}^{1+\alpha}+\Vert (\phi_2)' \Vert_{L_\infty}
\Vert h' \Vert_{C^\alpha})|x-y|^\alpha.
\end{align}
Moreover, $|(\phi_2 \circ h)'(x)|$ is bounded by $\Vert (\phi_2)' \Vert_{L_\infty} \Vert h' \Vert_{L_\infty}$.
As $\Vert h' \Vert_{L_\infty} \geq 1$,
these estimates imply that the composition operator $P_h$ is bounded on $C^{1+\alpha}(\mathbb S)$ 
with operator norm at most $\Vert h \Vert_{C^{1+\alpha}}^{1+\alpha}$.

By these bounds for $P_h$ on $\dot B^{1-\alpha}_{\infty,\infty}(\mathbb S)$ and $\ddot B^{1+\alpha}_{\infty,\infty}(\mathbb S)$
with equivalent norms, Proposition \ref{interpolation} implies that $P_h$ is a bounded operator acting on the real interpolation 
space in \eqref{represent} with operator norm bounded by
$c_\alpha \Vert h \Vert_{C^{1+\alpha}}$ for some constant $c_\alpha>0$.
Since $(P_h)^{-1}=P_{h^{-1}}$ and $h^{-1}$ also belongs to ${\rm Diff}^{1+\alpha}(\mathbb S)$,
we have that $(P_h)^{-1}$ is bounded on $C^Z(\mathbb S)$. 
Hence, $P_h$ is a Banach automorphism of $C^Z(\mathbb S)$.
\end{proof}

\begin{proof}[Proof of Claim \ref{claim1}]
It follows from Lemma \ref{composition} and \eqref{log} that $\log (h_\mu)'$ belongs to $C^Z(\mathbb S)$. 
This implies that $h_\mu=H(\mu)|_{\mathbb S}$ belongs to ${\rm Diff}^{1+Z}(\mathbb S)$.
\end{proof}

Finally, as an application of Lemma \ref{composition}, we confirm that 
${\rm Diff}^{1+Z}(\mathbb S)$ is a group.

\begin{proposition}\label{group}
${\rm Diff}^{1+Z}(\mathbb S)$ is a subgroup of ${\rm Diff}^{1+\alpha}(\mathbb S) \subset {\rm QS}(\mathbb S)$.
\end{proposition}

\begin{proof}
For $h_1, h_2 \in {\rm Diff}^{1+Z}(\mathbb S)$, we consider
$$
\log (h_1 \circ h_2)'=\log (h_1)' \circ h_2+\log (h_2)'.
$$
Since $\log (h_1)', \log (h_2)' \in C^Z(\mathbb S)$ and $\log (h_1)' \circ h_2 \in C^Z(\mathbb S)$
by Lemma \ref{composition}, we have $\log (h_1 \circ h_2)' \in C^Z(\mathbb S)$. 
For $h \in {\rm Diff}^{1+Z}(\mathbb S)$, we consider
$$
\log (h^{-1})'=-\log h' \circ h^{-1}.
$$
Again by Lemma \ref{composition}, we have $\log (h^{-1})' \in C^Z(\mathbb S)$.
\end{proof}

\begin{remark}
We equip ${\rm Diff}^{1+Z}(\mathbb S)$ with a right-invariant topology induced by the $C^Z$-norm. 
Namely, $h_n \in {\rm Diff}^{1+Z}(\mathbb S)$
converges to $h$ as $n \to \infty$ if $h_n \circ h^{-1} \to \rm id$ and $(h_n \circ h^{-1})' \to 1$ uniformly and moreover
$\log (h_n \circ h^{-1})' \to 0$ in $\Vert \cdot \Vert_{C^Z}$. Then, we expect that ${\rm Diff}^{1+Z}(\mathbb S)$ is
a {\it topological group}. In fact, Lemma \ref{composition} with the estimate of the operator norm of $P_h$ implies that 
${\rm Diff}^{1+Z}(\mathbb S)$ is a partial topological group in the sense that the group operations are continuous at
the identity (see \cite[Definition 1.2]{GS}).
To prove that ${\rm Diff}^{1+Z}(\mathbb S)$ is a topological group by
generalizing the continuity at $\rm id$ to any $h \in {\rm Diff}^{1+Z}(\mathbb S)$, it suffices to show that
the adjoint
$h \circ g \circ h^{-1}$ converges to $\rm id$ as $g  \in {\rm Diff}^{1+Z}(\mathbb S)$ converges to $\rm id$.
The same problem can be asked for the Teichm\"uller space $T^Z \cong  {\rm Diff}^{1+Z}_*(\mathbb S)$.
In this setting, a standard estimate of a Beltrami coefficient under composition might be useful.
\end{remark}

\section{The correspondence of $C^Z$ and $B^Z$}\label{4}

The following theorem is due to Zygmund. The boundedness of the operator can be
seen from the proof of \cite[Theorem 5.3]{Du}. A proof in a more general setting can be found in \cite[V. Proposition 8]{St}.

\begin{theorem}\label{zygmund}
A holomorphic function $\Phi$ in $B^Z(\mathbb D)$ extends continuously to $\mathbb S$, thereby defining a function 
$\phi=\Phi|_{\mathbb S}$ in $C^Z(\mathbb S)$.
The boundary extension operator $E:B^Z(\mathbb D) \to C^Z(\mathbb S)$ given in this way is
a Banach isomorphism onto its image.
\end{theorem}

The inverse $E^{-1}:E(B^Z(\mathbb D)) \to B^Z(\mathbb D)$ is given by the Poisson integral of $\phi$.
Moreover, this operator extends to all of $C^Z(\mathbb S)$ as the {\it Szeg\"o projection} defined by the Cauchy integral
$$
{\mathcal S}(\phi)(z)=\frac{1}{2\pi i} \int_{\mathbb S} \frac{\phi(\zeta)}{\zeta-z}\,d \zeta \quad (z \in \mathbb D)
$$
for $\phi \in C^Z(\mathbb S)$.
Then ${\mathcal S}(\phi)$ belongs to $B^Z(\mathbb D)$. 

On the other hand, 
the {\it Hilbert transform} is defined by the singular integral
$$
{\mathcal H}(\phi)(x)=\frac{1}{\pi i}\,{\rm p.v.} \int_{\mathbb S} \frac{\phi(\zeta)}{\zeta-x}\,d \zeta \quad (x \in \mathbb S).
$$
It is known that ${\mathcal H}$ maps $C^Z(\mathbb S)$ to itself.
This is a bounded linear operator on the Besov space $C^Z(\mathbb S)=\dot B^1_{\infty,\infty}(\mathbb S)$
(see \cite[Proposition 4.7]{GP} in this setting). Combined with Theorem \ref{zygmund}, this yields:

\begin{proposition}
The following hold:
\begin{itemize}
\item[$(1)$]
${\mathcal H}:C^Z(\mathbb S) \to C^Z(\mathbb S)$ is a Banach automorphism with ${\mathcal H} \circ {\mathcal H}=I$;
\item[$(2)$]
${\mathcal S}:C^Z(\mathbb S) \to B^Z(\mathbb D)$ is a bounded linear operator such that 
$E \circ {\mathcal S}:C^Z(\mathbb S) \to C^Z(\mathbb S)$ is a bounded projection onto 
$E(B^Z(\mathbb D))$ satisfying 
$E \circ {\mathcal S}=\frac{1}{2}(I +\mathcal H)$.
\end{itemize}
\end{proposition}
\smallskip

We also consider the boundary extension operator $E:B^Z(\mathbb D^*) \to C^Z(\mathbb S)$
on $\mathbb D^*$. Here, the corresponding space of holomorphic functions on $\mathbb D^*$ can be given simply by reflection as
$$
B^Z(\mathbb D^*)=\{\Psi \in {\rm Hol}(\mathbb D^*) \mid \Psi(z)=\Phi(z^*)^*,\ \Phi \in B^Z(\mathbb D)\}
$$ 
with norm $\Vert \Psi \Vert_{B^Z}=\Vert \Phi \Vert_{B^Z}$, where $z^*=1/\bar z$ is the reflection point of $z$
with respect to $\mathbb S$.
By defining the Szeg\"o projection
$$
{\mathcal S}^*(\phi)(z)=\frac{1}{2\pi i} \int_{\mathbb S} \frac{\phi(\zeta)}{\zeta-z}\,d \zeta \quad (z \in \mathbb D^*)
$$
for $\phi \in C^Z(\mathbb S)$, where the orientation of the line integral on $\mathbb S$ is taken counterclockwise, we see that
$E \circ {\mathcal S}^*$ coincides with the bounded projection $\frac{1}{2}(I -\mathcal H)$ 
onto $E(B^Z(\mathbb D^*))$.
Then the identifications $E(B^Z(\mathbb D)) \cong B^Z(\mathbb D)$ and $E(B^Z(\mathbb D^*)) \cong B^Z(\mathbb D^*)$
under both boundary extension operators $E$ yield
the topological direct sum decomposition
\begin{equation}\label{identify}
C^Z(\mathbb S) \cong B^Z(\mathbb D) \oplus B^Z(\mathbb D^*).
\end{equation}

For $(\mu_1, \mu_2) \in M^Z(\mathbb D) \times M^Z(\mathbb D^*)$, 
let $G(\mu_1, \mu_2): \widehat{\mathbb C} \to \widehat{\mathbb C}$ be the normalized quasiconformal self-homeomorphism with 
$\mu_G|_{\mathbb D} = \mu_1$, $\mu_G|_{\mathbb D^*} = \mu_2$.
The normalization is imposed by fixing $0$, $1$, and $\infty$; this ensures that the image of $\mathbb S$
under $G(\mu_1, \mu_2)$ is bounded. 

\begin{proposition}\label{simal}
If $\mu_1 \in M^Z(\mathbb D)$ and $\mu_2 \in M^Z(\mathbb D^*)$, then
$\log \gamma' \in C^Z(\mathbb S)$ for $\gamma=G(\mu_1, \mu_2)|_{\mathbb S}$. 
\end{proposition}

\begin{proof}
Let $\nu=\mu_2 \ast (\mu_1^*)^{-1}$, where $\mu^*$ denotes the reflection of a Beltrami coefficient 
$\mu$ with respect to $\mathbb S$.
Then $\nu \in M^Z(\mathbb D^*)$ by \cite[Proposition 4]{M3}, and we have
$G(\mu_1, \mu_2)= F_{\nu} \circ H(\mu_1)$. For $f=F_{\nu}|_{\mathbb S}$, $\log f'$ belongs to $C^Z(\mathbb S)$ by
Theorems \ref{equivalence} and \ref{zygmund}, and for $h=H(\mu_1)|_{\mathbb S}$, $\log h'$ belongs to $C^Z(\mathbb S)$ by Claim \ref{claim1}.
Then
$$
\log \gamma'=\log f' \circ h +\log h',
$$
and since $h \in {\rm Diff}^{1+Z}(\mathbb S)$, Lemma \ref{composition} shows that $\log \gamma' \in C^Z(\mathbb S)$.
\end{proof}

We define a map
\begin{equation}\label{Scase}
\widetilde \Lambda:M^Z(\mathbb D) \times M^Z(\mathbb D^*) \to C^Z(\mathbb S)
\end{equation}
by the correspondence $(\mu_1, \mu_2) \mapsto \log \gamma'$
for $\gamma=G(\mu_1, \mu_2)|_{\mathbb S}$.

\begin{lemma}\label{holo}
$\widetilde \Lambda$ is holomorphic.
\end{lemma}

For the proof of this lemma and also for later arguments, we consider the right translation in
$M^Z(\mathbb D)$ and $M^Z(\mathbb D^*)$ with respect to the group structure.
For any $\nu \in M^Z(\mathbb D)$, define 
$r_{\nu}:M^Z(\mathbb D) \to M^Z(\mathbb D)$ by $\mu \mapsto \mu \ast \nu$, where $\mu \ast \nu$ denotes
the complex dilatation of $H(\mu) \circ H(\nu)$.
This is a biholomorphic automorphism of $M^Z(\mathbb D)$ (see \cite[Lemma 5]{M3}).
The same is true on $M^Z(\mathbb D^*)$.
Then, by the skew-diagonal action 
$$
r_{\nu}(\mu_1, \mu_2) = (\mu_1 \ast \nu, \mu_2 \ast \nu^*),
$$
it also acts biholomorphically on $M^Z(\mathbb D) \times M^Z(\mathbb D^*)$.

From the facts that the Schwarzian derivative map $S:M^Z(\mathbb D^*) \to A^Z(\mathbb D)$ is a holomorphic split submersion
onto its image $S(M^Z(\mathbb D^*))$
and the Bers embedding $\alpha:T^Z \to A^Z(\mathbb D)$ is a biholomorphic homeo\-morphism onto $S(M^Z(\mathbb D^*))$
 (see \cite[Theorem 3]{M3}),
we have a local holomorphic right inverse to the Teichm\"uller projection 
$\pi:M^Z(\mathbb D^*) \to T^Z=M^Z(\mathbb D^*)/\sim$.
Then we can project the right translation $r_{\nu}$ of $M^Z(\mathbb D^*)$ down to a right translation on $T^Z$
as a biholomorphic automorphism $R_{[\nu]}$ induced by $[\nu] \in T^Z$, namely
$R_{[\nu]}:[\mu] \mapsto [\mu]\ast[\nu]:=[\mu \ast \nu]$ (see \cite[Remark 1]{M3}).
Moreover, this extends to the biholomorphic automorphism $R_{[\nu]}$ of 
the product of the Teichm\"uller spaces $(M^Z(\mathbb D)/\sim) \times (M^Z(\mathbb D^*)/\sim)$
defined by 
\begin{equation}\label{product}
R_{[\nu]}([\mu_1], [\mu_2]) = ([\mu_1] \ast [\nu], [\mu_2] \ast [\nu^*]).
\end{equation}

For $h \in {\rm Diff}^{1+Z}(\mathbb S)$, define an affine translation $Q_h$ on $C^Z(\mathbb S)$ 
by $Q_h(\phi) = P_h(\phi) + \log h'$, where the composition operator $P_h$ 
is a Banach automorphism by Lemma \ref{composition}.
The right translation $r_\nu$ on $M^Z(\mathbb D) \times M^Z(\mathbb D^*)$ and the affine translation $Q_h$ on $C^Z(\mathbb S)$ 
satisfy the following relation under $\widetilde \Lambda$. The argument is the same as that for \cite[Proposition 5.1]{WM-3}.

\begin{proposition}\label{translation}
We have
$$
\widetilde \Lambda \circ r_\nu = Q_h \circ \widetilde \Lambda
$$
for $h=H(\nu)|_{\mathbb S}$ with $\nu \in M^Z(\mathbb D)$.
\end{proposition}

\begin{proof}[Proof of Lemma \ref{holo}]
By the Hartogs theorem for Banach spaces (see \cite[\S 14.27]{Ch}), to see that $\widetilde \Lambda$ is holomorphic
it suffices to show that $\widetilde \Lambda$ is separately holomorphic.
Namely, fix $\nu \in M^Z(\mathbb D)$ and prove that $\widetilde \Lambda(\nu,\,\cdot \,)$ is holomorphic. 
The other case is treated in the same way.

Let $h=H(\nu)|_{\mathbb S} \in {\rm Diff}^{1+Z}(\mathbb S)$.
For the affine translation $Q_{h}$ on $C^Z(\mathbb S)$ induced by $h$,
Proposition \ref{translation} gives 
$\widetilde \Lambda \circ r_{\nu} = Q_{h} \circ \widetilde \Lambda$.
This relation yields
the useful representation
\begin{equation}
\widetilde \Lambda(\nu,\,\cdot \,)=Q_{h} \circ \widetilde \Lambda\!\left(0,r_{(\nu^*)^{-1}}(\,\cdot \,)\right).
\end{equation}
Here, 
$\widetilde \Lambda(0,\,\cdot \,)$ can be regarded as 
the pre-Schwarzian derivative map $L:M^Z(\mathbb D^*) \to B^Z(\mathbb D)$ defined by $\mu \mapsto \log (F_\mu)'$
composed with the boundary extension $E$, that is,
\begin{equation}\label{pre-Schwarz}
\widetilde \Lambda(0,\mu)=E(L(\mu)) \quad (\mu \in M^Z(\mathbb D^*)).
\end{equation}
Since $E$ is a bounded linear operator, and $L$ and $r_{(\nu^*)^{-1}}$ are holomorphic, we conclude that
$\widetilde \Lambda(\nu,\,\cdot \,)$ is holomorphic.
\end{proof}

\section{Solution to Claim \ref{claim2}}\label{5}

In this section, we address Claim \ref{claim2}.
We give a proof of this claim 
using simultaneous uniformization in the theory of absolutely continuous Teichm\"uller spaces (see \cite{M5}).
This argument involves function spaces defined on the real line $\mathbb R$, into which the Teichm\"uller spaces defined on
the upper and lower half-planes $\mathbb H$ and $\mathbb H^*$ are embedded.
To this end, 
we lift functions defined on $\mathbb D \setminus \{0\}$, $\mathbb D^* \setminus \{\infty\}$, and $\mathbb S$ to 
periodic functions on $\mathbb H$, $\mathbb H^*$, and $\mathbb R$, respectively, by the universal cover. See diagram \eqref{univ}.

We set the following spaces. Here, $M(\mathbb H^*)$ is the space of Beltrami coefficients on $\mathbb H^*$, and
$B(\mathbb H)$ is the space of Bloch functions $\Phi$ on $\mathbb H$,
consisting of functions holomorphic on $\mathbb H$ with $\sup_{{\rm Im}\, z > 0}
({\rm Im}\, z)|\Phi'(z)|<\infty$.
\begin{align*}
M^Z(\mathbb H^*)&=\{ \mu \in M(\mathbb H^*) \mid \Vert \mu \Vert_Z 
= \underset{{\rm Im}\, z < 0}{\mathrm{ess\,sup}}\, (|{\rm Im}\, z|^{-1} \lor 1)|\mu(z)| < \infty \}, \\
&M^Z_{\rm per}(\mathbb H^*) = \{\mu \in M^Z(\mathbb H^*) \mid \mu(z+2\pi) = \mu(z)\ (\forall z \in \mathbb H^*)\}; \\
C^Z(\mathbb R) &= \{\phi \in C(\mathbb R) \mid \Vert \phi \Vert_{C^Z}= \\
&\qquad \qquad \qquad \sup_{x \in \mathbb R,\,t > 0} (2t)^{-1}|\phi(x + t) + \phi(x - t) - 2\phi(x)| < \infty\}, \\
&C^Z_{\rm per}(\mathbb R) = \{\phi \in C^Z(\mathbb R) \mid \phi(x+2\pi) = \phi(x)\ (\forall x \in \mathbb R)\};\\
B^Z(\mathbb H)&=\{ \Phi \in B(\mathbb H) \mid \Vert \Phi \Vert_{B^Z}=\sup_{{\rm Im}\, z > 0}
({\rm Im}\, z)|\Phi''(z)|+|\Phi'(i)|<\infty \},\\
&B^Z_{\rm per}(\mathbb H) =\{ \Phi \in B^Z(\mathbb H) \mid \Phi(z+2\pi) = \Phi(z)\ (\forall z \in \mathbb H)\}.
\end{align*}
Moreover, $M^Z_{\rm per}(\mathbb H)$ and $B^Z_{\rm per}(\mathbb H^*)$ are defined analogously.
Then, we have the following identifications:
\begin{align*}
&M^Z_{\rm per}(\mathbb H^*) \cong M^Z(\mathbb D^*), \quad M^Z_{\rm per}(\mathbb H) \cong M^Z(\mathbb D); \quad
C^Z_{\rm per}(\mathbb R) \cong C^Z(\mathbb S);\\
&B^Z_{\rm per}(\mathbb H) \cong B^Z(\mathbb D), \quad B^Z_{\rm per}(\mathbb H^*) \cong B^Z(\mathbb D^*).
\end{align*}

The Teichm\"uller space $T^Z$ defined by $M^Z(\mathbb D^*)/\sim$ can also be given by
\linebreak
$M^Z_{\rm per}(\mathbb H^*)/\sim$ under the analogously defined Teichm\"uller equivalence $\sim$.
The topological direct sum decomposition \eqref{identify} yields
\begin{equation}\label{decompose}
C^Z_{\rm per}(\mathbb R) \cong B^Z_{\rm per}(\mathbb H) \oplus B^Z_{\rm per}(\mathbb H^*).
\end{equation}
The real Banach subspace of $C^Z_{\rm per}(\mathbb R)$ consisting of all real-valued functions is denoted by
${\rm Re}\, C^Z_{\rm per}(\mathbb R)$.

For $(\mu_1, \mu_2) \in M^Z_{\rm per}(\mathbb H) \times M^Z_{\rm per}(\mathbb H^*)$, 
let $G(\mu_1, \mu_2): \mathbb C \to \mathbb C$ be the normalized quasiconformal self-homeomorphism of $\mathbb C$ with 
$\mu_G|_{\mathbb H} = \mu_1$, $\mu_G|_{\mathbb H^*} = \mu_2$.
The normalization is imposed by fixing $0$, $2\pi$, and $\infty$.
Since the Beltrami coefficients are $2\pi$-periodic, this normalization implies
$$
G(\mu_1,\mu_2)(z+2\pi)= G(\mu_1,\mu_2)(z)+2\pi.
$$
A similar map was defined before on $\mathbb D$ and $\mathbb D^*$ with the same notation.
However, in the present case,
the boundary value $\gamma = G(\mu_1, \mu_2)|_{\mathbb R}$ is determined by the pair of
Teichm\"uller equivalence classes $([\mu_1], [\mu_2])$ (see \cite[Proposition 4.1]{WM-3}), 
and $\log \gamma'$ belongs to $C^Z_{\rm per}(\mathbb R)$ by Proposition \ref{simal}. 
Thus, the map
$$
\Lambda: T^Z \times T^Z_* \to C^Z_{\rm per}(\mathbb R)
$$
for $T^Z=M^Z_{\rm per}(\mathbb H)/\sim$ and $T^Z_*=M^Z_{\rm per}(\mathbb H^*)/\sim$ is induced in this way.

By applying the local holomorphic right inverse to the Teichm\"uller projection, 
we see from Lemma \ref{holo} that $\Lambda$ is holomorphic. Moreover,
the biholomorphic automorphism $R_{[\nu]}$ of $T^Z \times T^Z_*$ is defined for $[\nu] \in T^Z$ in the same way as \eqref{product},
and Proposition \ref{translation} translates into 
\begin{equation}\label{relation}
\Lambda \circ R_{[\nu]} = Q_h \circ \Lambda
\end{equation}
for $h=G(\nu,\bar \nu)|_{\mathbb R}=H(\nu)|_{\mathbb R}$.
Here, the reflection of a Beltrami coefficient $\nu$ with respect to $\mathbb R$ is denoted by $\bar \nu$, and 
the quasiconformal self-homeomorphism of $\mathbb H$ with the complex dilatation $\nu$ and
the normalization fixing $0$, $2\pi$, and $\infty$ is denoted by $H(\nu)$.

\begin{proposition}\label{prop.5.1}
$\Lambda: T^Z \times T^Z_* \to C^Z_{\rm per}(\mathbb R)$ is a holomorphic injection satisfying
$\Lambda \circ R_{[\nu]} = Q_h \circ \Lambda$ for $h=H(\nu)|_{\mathbb R}$ with $\nu \in M^Z_{\rm per}(\mathbb H)$.
\end{proposition}

\begin{proof}
Only the injectivity of $\Lambda$ remains to be proved.
Suppose that 
\linebreak 
$\Lambda([\mu_1],[\nu_1])=\Lambda([\mu_2],[\nu_2])$. Then,
$G(\mu_1,\nu_1)|_{\mathbb R}=G(\mu_2,\nu_2)|_{\mathbb R}$ by the normalization fixing the three points on $\mathbb R$.
This implies that $[\mu_1]=[\mu_2]$ and $[\nu_1]=[\nu_2]$, which can be verified
by the same proof as \cite[Proposition 4.1]{WM-3}.
Hence, $\Lambda$ is injective.
\end{proof}

Let $\Delta(T^Z)=\{([\mu],[\bar \mu]) \in T^Z \times T^Z_* \mid [\mu] \in T^Z\}$ be the anti-diagonal axis of $T^Z \times T^Z_*$,
which is a real-analytic submanifold real-analytically equivalent to $T^Z$. 

Now, Claim \ref{claim2} is deduced from the following theorem. 

\begin{theorem}\label{real}
The restriction of $\Lambda$ to the anti-diagonal satisfies
$$
 \Lambda(\Delta(T^Z))
 = {\rm Re}\,C^Z_{\rm per}(\mathbb R).
$$
In particular, the image of $\Lambda$ contains
${\rm Re}\,C^Z_{\rm per}(\mathbb R)$.
\end{theorem}

\begin{proof}[Proof of Claim \ref{claim2}]
We lift a quasisymmetric homeomorphism in ${\rm Diff}^{1+Z}(\mathbb S)$ to $h:\mathbb R \to \mathbb R$
such that $\log h' \in {\rm Re}\,C^Z_{\rm per}(\mathbb R)$. By Theorem \ref{real}, this lies in the image of $\Delta(T^Z)$
under $\Lambda$,
which implies that there exists $\mu \in M^Z_{\rm per}(\mathbb H)$ such that $\Lambda([\mu],[\bar \mu])=\log h'$.
Hence, $h=G(\mu,\bar \mu)|_{\mathbb R}=H(\mu)|_{\mathbb R}$, and thus the claim follows.
\end{proof}

Theorem \ref{real} is proved by applying the following lemma.

\begin{lemma}\label{origin}
$\Lambda$ is locally biholomorphic at the origin $([0],[0]) \in T^Z \times T^Z_*$.
\end{lemma}

\begin{proof}
We show that the derivative of $\Lambda$ at $([0],[0])$ is a surjective isomorphism.
Then the inverse mapping theorem (see \cite[\S 7.18]{Ch}) yields the statement.
If we fix the first coordinate as $[0]$, then taking the Teichm\"uller projection for \eqref{pre-Schwarz}, we see that
$\Lambda([0],\cdot)$ is nothing but
the pre-Bers embedding $\beta:T^Z_* \to B^Z_{\rm per}(\mathbb H)$, where $B^Z_{\rm per}(\mathbb H)$ is identified
with its image in $C^Z_{\rm per}(\mathbb R)$ under the boundary extension isomorphism $E$. 
We note that in the case of the half-plane, the pre-Schwarzian derivative map 
$L:M^Z_{\rm per}(\mathbb H^*) \to B^Z_{\rm per}(\mathbb H)$ is factored through the Teichm\"uller projection to $\beta$,
as with the Bers embedding. See \cite[Section 7]{WM-4}.

It follows that
the derivative $d_{([0],[0])}\Lambda$ maps the tangent subspace along the second coordinate isomorphically 
onto $B^Z_{\rm per}(\mathbb H)$.
In the same way, $d_{([0],[0])}\Lambda$ maps the tangent subspace along the first coordinate isomorphically
onto $B^Z_{\rm per}(\mathbb H^*)$. Since $C^Z_{\rm per}(\mathbb R)=B^Z_{\rm per}(\mathbb H) \oplus B^Z_{\rm per}(\mathbb H^*)$
by \eqref{decompose}, we have that $d_{([0],[0])}\Lambda$ is a surjective isomorphism.
\end{proof}

\begin{remark}
We have transferred our arguments from $\mathbb D$ and $\mathbb S$ to $\mathbb H$ and $\mathbb R$ with periodicity.
The reason for doing so is that we have the following advantages:
(1) The pre-Schwarzian derivative map $L$, or more generally, $\Lambda$ induced by $\widetilde \Lambda$ in \eqref{Scase},
is well defined as a map from the Teichm\"uller space; 
(2) For a suitably given function $\phi$ on $\mathbb R$,
a homeomorphism of $\mathbb R$ can be constructed simply by $\int_0^x \exp(\phi(t))\,dt$.
\end{remark}

\begin{proof}[Proof of Theorem \ref{real}]
Consider the involution
$$
\iota([\mu_1],[\mu_2])=([\bar\mu_2],[\bar\mu_1])
$$
on $T^Z\times T^Z_*$. Then, $\Delta(T^Z)$ is the fixed point locus of $\iota$.
By reflection with respect to $\mathbb R$, we have
$\Lambda\circ\iota=\overline{\Lambda}$.
Since $\phi \in {\rm Re}\,C^Z_{\rm per}(\mathbb R)$ if $\bar \phi=\phi$ for $\phi \in C^Z_{\rm per}(\mathbb R)$ and since
$\Lambda$ is injective by Proposition 5.1, it follows that
\begin{equation}\label{reallocus}
 \Lambda^{-1}({\rm Re}\,C^Z_{\rm per}(\mathbb R))
 =\operatorname{Fix}(\iota)
 =\Delta(T^Z).
\end{equation}

By Lemma \ref{origin}, there is a neighborhood $U$ of
$([0],[0])$ in $T^Z\times T^Z_*$ such that $\Lambda$ maps $U$
biholomorphically onto a neighborhood
$V=\Lambda(U)$ of the origin in $C^Z_{\rm per}(\mathbb R)$.
Set $A=\Lambda(\Delta(T^Z))$.
Then
\[
 A_0=V\cap{\rm Re}\,C^Z_{\rm per}(\mathbb R)
\]
is an open neighborhood of the origin in
${\rm Re}\,C^Z_{\rm per}(\mathbb R)$ and, by \eqref{reallocus},
$A_0\subset A$.

We equip ${\rm Re}\,C^Z_{\rm per}(\mathbb R)$ with the group
structure transported from the group of equivariant
$C^{1+Z}$-diffeomorphisms of $\mathbb R$. More explicitly,
for $\phi\in{\rm Re}\,C^Z_{\rm per}(\mathbb R)$, we choose its
representative, modulo constants, so that
$\int_0^{2\pi}e^{\phi(t)}\,dt=2\pi$,
and put
\[
 h_\phi(x)=\int_0^x e^{\phi(t)}\,dt.
\]
Then $h_\phi(x+2\pi)=h_\phi(x)+2\pi$ and
$\log h_\phi'=\phi$. For $\psi,\phi\in
{\rm Re}\,C^Z_{\rm per}(\mathbb R)$, the group operation is
\[
 \psi*\phi
 =Q_{h_\phi}(\psi)
 =\psi\circ h_\phi+\log h_\phi'
 =\log (h_\psi\circ h_\phi)'.
\]

Under the identification $\Delta(T^Z)\cong T^Z$, the group
operation is given by
\[
 ([\mu],[\bar\mu])*([\nu],[\bar\nu])
 =
 R_{[\nu]}([\mu],[\bar\mu]).
\]
Hence, by \eqref{relation}, the restriction
$\Lambda:\Delta(T^Z)\to
{\rm Re}\,C^Z_{\rm per}(\mathbb R)$
is a group homomorphism. Therefore $A$ is a subgroup of
${\rm Re}\,C^Z_{\rm per}(\mathbb R)$.

Since $A$ contains the open neighborhood $A_0$ of the identity,
$A$ is open. Indeed, every point $\phi \in A$ has an open neighborhood
obtained from $A_0$ by a right translation $Q_{h_\phi}$.
Each $Q_h$ for $\log h' \in {\rm Re}\,C^Z_{\rm per}(\mathbb R)$ is an affine homeomorphism by Lemma \ref{composition},
and hence every right coset of $A$ in ${\rm Re}\,C^Z_{\rm per}(\mathbb R)$ is open. Consequently,
the complement of $A$, being the union of all the other right
cosets, is also open. Thus $A$ is both open and closed in
${\rm Re}\,C^Z_{\rm per}(\mathbb R)$.

Since ${\rm Re}\,C^Z_{\rm per}(\mathbb R)$ is a real Banach
space and hence connected, we conclude that
$A={\rm Re}\,C^Z_{\rm per}(\mathbb R)$.
\end{proof}

\section{The real-analytic structure of $T^Z$}\label{6}

We have shown in Proposition \ref{prop.5.1} that $\Lambda$ is a holomorphic injection.
In this section, we assert more:

\begin{theorem}\label{biholo}
The map $\Lambda: T^Z \times T^Z_* \to C^Z_{\rm per}(\mathbb R)$ 
is a biholomorphism onto its image.
\end{theorem}

By Theorem \ref{real}, 
$\Lambda(\Delta(T^Z))={\rm Re}\,C^Z_{\rm per}(\mathbb R)$.
Moreover, $\Lambda$ is locally biholomorphic on a neighborhood $U$ of $([0],[0]) \in T^Z \times T^Z_*$ 
by Lemma \ref{origin}, and, by right translations and affine translations satisfying \eqref{relation},  
this extends to a neighborhood of $\Delta(T^Z)$.
Thus, $\Lambda$ is a biholomorphism from a neighborhood of $\Delta(T^Z)$ onto its image, which is a neighborhood of 
${\rm Re}\,C^Z_{\rm per}(\mathbb R)$ in $C^Z_{\rm per}(\mathbb R)$.
This already yields the following corollary to Theorem \ref{biholo}
without establishing the biholomorphy of $\Lambda$ on the entire $T^Z \times T^Z_*$.

\begin{corollary}\label{real-analytic}
The Teichm\"uller space $T^Z$ is real-analytically equivalent to 
the real Banach space ${\rm Re}\,C^Z_{\rm per}(\mathbb R)$.
\end{corollary}

\begin{corollary}\label{contractible}
$T^Z$ is real-analytically contractible; that is,
$T^Z$ contracts to a point via a real-analytic homotopy.
\end{corollary}

\begin{remark}
Our arguments also apply to the Teichm\"uller space $T^\alpha \cong$
\linebreak
${\rm Diff}^{1+\alpha}_*(\mathbb S)$ 
for $0 <\alpha <1$ investigated in \cite{M1}.
Hence, $T^\alpha$ is real-analytically equivalent to 
the real Banach space ${\rm Re}\,C^\alpha_{\rm per}(\mathbb R)$, and
$T^\alpha$ is real-analytically contractible. Topological contractibility of $T^\alpha$ was
proved by using a global continuous section to the Teichm\"uller projection constructed by
the barycentric extension of the elements of ${\rm Diff}^{1+\alpha}(\mathbb S)$ (see \cite[Theorem 1.1]{M0}).
\end{remark}

\begin{proof}[Proof of Theorem \ref{biholo}]
We outline the proof, since the full arguments appear in \cite[Theorem 7.1]{M5} and \cite[Theorem 9.1]{M6} for different
Teichm\"uller spaces, and essentially the same method applies here.

Since $\Lambda$ is a holomorphic injection by Proposition \ref{prop.5.1},
it suffices to prove that the inverse $\Lambda^{-1}$ is holomorphic.
By the inverse mapping theorem as before, it is enough to show that the derivative $d\Lambda$ is surjective
at any point.
The proof of the surjectivity of
$$
d_{([\mu],[\nu])}\Lambda: {\mathscr T}_{([\mu],[\nu])}(T^Z \times T^Z_*) \to C^Z_{\rm per}(\mathbb R)
$$
on the tangent space ${\mathscr T}_{([\mu],[\nu])}(T^Z \times T^Z_*)$ at $([\mu],[\nu]) \in T^Z \times T^Z_*$
proceeds in the following steps.

(i)
The image of the tangent space
$$
{\mathscr T}_{([\mu],[\nu])}(T^Z \times T^Z_*) \cong B^Z_{\rm per}(\mathbb H^*) \oplus B^Z_{\rm per}(\mathbb H)
$$
under $d_{([\mu],[\nu])}\Lambda$ is the algebraic sum
$
P_{h_\mu}B^Z_{\rm per}(\mathbb H)+ P_{h_\nu}B^Z_{\rm per}(\mathbb H^*),
$
where $h_\mu=H(\mu)|_{\mathbb R}$ and $h_\nu=H(\nu)|_{\mathbb R}$.

(ii)
If $\Lambda([\mu_0],[\nu_0])$ lies in the real subspace $i\,{\rm Re}\,C^Z_{\rm per}(\mathbb R)$ of $C^Z_{\rm per}(\mathbb R)$, then
\linebreak 
$d_{([\mu_0],[\nu_0])}\Lambda$ is surjective onto $i\,{\rm Re}\,C^Z_{\rm per}(\mathbb R)$.
It follows from this and (i) that 
\linebreak
$d_{([\mu_0],[\nu_0])}\Lambda$ is surjective for such a pair $([\mu_0],[\nu_0])$.

(iii) 
We establish (ii)
by relying on the result in \cite[Theorem 6.1]{WM-3} concerning the biholomorphy of $\Lambda$ on a certain domain
in the product of the BMO Teichm\"uller spaces. This domain contains $T^Z \times T^Z_*$, and
the surjectivity of $d\Lambda$ for this larger space can be invoked.
We also use the Kellogg--Warschawski theorem to see that the Riemann mapping parametrization has the required regularity.

(iv)
If $\log \gamma_0' \in i\,{\rm Re}\,C^Z_{\rm per}(\mathbb R)$, then $\gamma_0$ is parametrized by arc-length.
Any curve $\gamma$ with $\log \gamma' = \Lambda([\mu],[\nu])$ can be reparametrized from that of arc-length, and
this is achieved by some homeomorphism $h=H(\lambda)|_{\mathbb R}$ with $\lambda \in M_{\rm per}^Z(\mathbb H)$.

(v)
For $\Lambda([\mu_0],[\nu_0]) = \log \gamma_0'$,
we move $\log \gamma_0'$ by the affine translation $Q_{h}$, and
translate $([\mu_0],[\nu_0])$ in parallel by the right translation $R_{[\lambda]}$
under the correspondence \eqref{relation}.
Then by (ii), 
the surjectivity of $d_{([\mu],[\nu])}\Lambda$ at any point $([\mu],[\nu]) \in T^Z \times T^Z_*$ follows.
\end{proof}

\section{Another proof of $M^Z \Rightarrow A^Z$}\label{7}

In this section, as an appendix, we give a technical but constructive proof of the implication $(1) \Rightarrow (3)$ in
Theorem \ref{equivalence}. It also makes the estimate of the norm of a Schwarzian derivative
in terms of the norm of a Beltrami coefficient clear.
This method was used in \cite[Theorem 4.1]{M1} to prove $(1) \Rightarrow (2)$ in the case of $T^\alpha$,
but it can also be applied to the present case.
We state the results more generally for $0<\alpha<2$, though our interest is in $\alpha=1$.
For ease of presentation, we replace the reciprocal of the hyperbolic density $1-|\zeta|^2$ in $\mathbb D$ by $1-|\zeta|$; since
$1-|\zeta| \leq 1-|\zeta|^2 \leq 2(1-|\zeta|)$, this causes no difficulty.

\begin{theorem}\label{alpha2}
For every $\alpha \in (0,2)$, there exists a constant $C>0$ depending only on $\alpha$
such that if a Beltrami coefficient $\mu \in M(\mathbb D^*)$ satisfies
$$
\Vert \mu \Vert_{(\alpha)}:=\underset{|z|>1}{\mathrm{ess\,sup}}\, ((|z|-1)^{-\alpha} \lor 1)|\mu(z)|<\infty,
$$
then for the quasiconformal self-homeomorphism $F_\mu$ of $\widehat{\mathbb C}$ with complex dilatation $\mu$ on $\mathbb D^*$
and $0$ on $\mathbb D$, the Schwarzian derivative of $F_\mu|_{\mathbb D}$ satisfies
$$
(1-|\zeta|)^{2-\alpha}|S_{F_\mu|_{\mathbb D}}(\zeta)| \leq C \Vert \mu \Vert_{(\alpha)}
$$
for every $\zeta \in \mathbb D$.
\end{theorem}

We prepare the following two lemmas for the proof of this theorem.

\begin{lemma}\label{recurrencelemma}
For any $\alpha$ with $0< \alpha < 2$,
there exists a constant $\lambda$ $(0 < \lambda <1)$ depending only on $\alpha$ such that
if a sequence $\{s_n\}_{n=0}^\infty$ of positive real numbers satisfies the recurrence relation
$$
\left(\frac{1}{1+s_{n-1}} \right)^2 {s_n}^\alpha=\lambda^n 
$$
for every $n \geq 1$ and $s_0=1$, 
then $\{s_n\}$ is increasing and diverges to $+\infty$.
\end{lemma}

\begin{proof}
The recurrence relation is equivalent to
$$
s_n=\lambda^{\frac{n}{\alpha}}(1+s_{n-1})^{\frac{2}{\alpha}}
$$
for every $n \geq 1$ with $s_0=1$. Modifying this formula, consider
\begin{equation}\label{recurrence}
s'_n=\lambda^{\frac{n}{\alpha}}{s'_{n-1}}^{\frac{2}{\alpha}}
\end{equation}
for every $n \geq 2$, with the initial value $s'_1=s_1=(4\lambda)^{1/\alpha}$.
It is easy to see that $s_n \geq s'_n$ for every $n \geq 1$. Hence,
$\lim_{n \to \infty} s'_n=+\infty$ implies $\lim_{n \to \infty} s_n=+\infty$.
Moreover, if $\{s'_n\}$ is increasing, then so is $\{s_n\}$.

Set $b_n=s'_{n+1}/s'_n$. Then \eqref{recurrence} becomes
$$
b_n=\lambda^{\frac{1}{\alpha}} (b_{n-1})^{\frac{2}{\alpha}}
$$
for every $n \geq 2$, and 
$$
b_1=\frac{s'_2}{s'_1}=\frac{\lambda^{\frac{2}{\alpha}}(4\lambda)^{\frac{2}{\alpha^2}}}{(4\lambda)^{\frac{1}{\alpha}}}.
$$
Taking logarithms yields
\begin{align*}
\log b_n&=\frac{2}{\alpha} \log b_{n-1} +\frac{1}{\alpha} \log\lambda,\\
\log b_1&=\left(\frac{2}{\alpha^2}+\frac{1}{\alpha}\right) \log \lambda+
\left(\frac{2}{\alpha^2}-\frac{1}{\alpha}\right)\log 4.
\end{align*}
From the first recurrence equation, if
$$
\log b_1 > \frac{-\log \lambda}{2-\alpha},
$$
then $\log b_n$ is positive and uniformly bounded away from $0$ for all $n \geq 1$.
Choosing $\lambda<1$ sufficiently close to $1$ achieves this.
For instance, one may take 
$\lambda>(1/\sqrt{2})^{(2-\alpha)^2}$.
Thus $\{s'_n\}$ is increasing and diverges to $+\infty$.
\end{proof}

\begin{lemma}\label{decomposition}
For a finite sequence of real numbers
$$
1=R_{-1}<R_0<R_1< \cdots <R_{N}<R_{N+1}=+\infty,
$$
let $A_i=\{R_i < |z| < R_{i+1}\}$ be an annulus (possibly degenerate)
in $\mathbb D^*$ 
for each $i=-1,0,\ldots,N$.  
Define
$$
\mu_i(z)=\left\{
\begin{array}{l} \mu(z) \,\quad (z \in A_i),\\[2pt]
\ 0  \qquad\, (z \in \widehat{\mathbb C} \setminus A_i),
\end{array}
\right.
$$ 
for $\mu \in M(\mathbb D^*)$, and set $k_i=\Vert \mu_i \Vert_{L_\infty}$. Then
$$
|S_{F_\mu|_\mathbb D}(\zeta)| 
\leq 12 \sum_{i=-1}^{N} \frac{k_i}{(R_i-|\zeta|)^2} \qquad (\zeta \in \mathbb D).
$$
\end{lemma}

\begin{proof}
First take a quasiconformal self-homeomorphism $F_{N}=F_{\mu_{N}}$ of $\widehat{\mathbb C}$ with complex dilatation 
$\mu_{N}$ and consider the push-forward $\tilde \mu_{N-1}=(F_{N})_* \mu_{N-1}$.
Then take 
$F_{N-1}=F_{\tilde \mu_{N-1}}$ and the push-forward 
$\tilde \mu_{N-2}=(F_{N-1}\circ F_{N})_* \mu_{N-2}$,
and continue inductively. Namely, for each $i \geq 0$,
let $F_{i}=F_{\tilde \mu_{i}}$ be the quasiconformal self-homeomorphism with complex dilatation 
$\tilde \mu_{i}$, and let 
$\tilde \mu_{i-1}=(F_{i}\circ \cdots \circ F_{N})_* \mu_{i-1}$ be the push-forward 
of $\mu_{i-1}$ by $F_{i}\circ \cdots \circ F_{N}$.
Finally, choose $F_{-1}=F_{\tilde \mu_{-1}}$ so that 
$F_{-1} \circ \cdots \circ F_{N}=F_\mu$.

The chain rule for Schwarzians shows that
\begin{align}\label{chain}
&\quad S_{F_\mu|_{\mathbb D}}(\zeta)\\
&=S_{F_{N}}(\zeta)+ \cdots +
S_{F_{-1}}(F_{0} \circ \cdots \circ F_{N}(\zeta))(F_{0} \circ \cdots \circ F_{N})'(\zeta)^2\\
&=
\sum_{i=-1}^{N} S_{F_i}(F_{i+1} \circ \cdots \circ F_{N}(\zeta))(F_{i+1} \circ \cdots \circ F_{N})'(\zeta)^2.
\end{align}
By the standard estimate for the Schwarzian derivative,
the conformal homeomorphism $F_{N}$ of the disk 
$\Omega_{N}=\{|\zeta|<R_{N}\}$ into $\widehat{\mathbb C}$ satisfies 
\begin{equation}\label{i=N}
|S_{F_{N}}(\zeta)|
\leq \frac{6k_{N}R_{N}^2}{(R_{N}^2-|\zeta|^2)^2} \leq \frac{6k_{N}}{(R_{N}-|\zeta|)^2}
\end{equation}
for every $\zeta \in \Omega_{N}$.

On the other hand, the conformal homeomorphism $F_{i}$ of the 
quasidisk $\Omega_i$ with hyperbolic density $\rho_{\Omega_i}$ (of curvature $-4$) into $\mathbb C$ for $-1 \leq i \leq N-1$, 
where $\Omega_i$ is the image of the disk $\{|\zeta|<R_i\}$ under
$F_{i+1} \circ \cdots \circ F_{N}$, satisfies 
$$
|S_{F_i}(\omega)| \leq 2 \cdot 6k_i\, \rho_{\Omega_i}(\omega)^2
$$
for every $\omega \in \Omega_i$.
Setting $\omega=F_{i+1} \circ \cdots \circ F_{N}(\zeta)$, we obtain
\begin{align}\label{i=other}
&\quad \ |S_{F_i}(F_{i+1} \circ \cdots \circ F_{N}(\zeta))(F_{i+1} \circ \cdots \circ F_{N})'(\zeta)^2|\\
&\leq 12k_i \,\rho_{\Omega_i}(F_{i+1} \circ \cdots \circ F_{N}(\zeta))^2|(F_{i+1} 
\circ \cdots \circ F_{N})'(\zeta)|^2\\
&=\frac{12k_i R_i^2}{(R_i^2-|\zeta|^2)^2} \leq  \frac{12k_i}{(R_i-|\zeta|)^2}
\end{align}
for every $\zeta \in \mathbb D$.
Plugging \eqref{i=N} and \eqref{i=other} into \eqref{chain} gives
$$
|S_{F_\mu|_\mathbb D}(\zeta)| \leq 12 \sum_{i=-1}^{N} \frac{k_i}{(R_i-|\zeta|)^2}
$$
for every $\zeta \in \mathbb D$.
\end{proof}
\medskip

\begin{proof}[Proof of Theorem \ref{alpha2}]
We estimate $(1-|\zeta|)^2|S_{F_\mu|_\mathbb D}(\zeta)|$ for each fixed $\zeta \in \mathbb D$.
Let $\tau=1-|\zeta|$ and $\ell=\Vert \mu \Vert_{(\alpha)}$.
If $\ell=0$, the assertion is trivial. Hence, we assume that $\ell>0$.
With the initial condition $t_0=\tau$, define $\{t_n\}_{n\geq 1}$ inductively by
\begin{equation}\label{t-rec}
\left(\frac{\tau}{\tau+t_{n-1}}\right)^2 \cdot \ell {t_n}^\alpha=
\lambda^n \cdot \ell \tau^\alpha 
\end{equation}
for some constant $\lambda$ with $0<\lambda<1$.
This is equivalent to
$$
\left(\frac{1}{1+s_{n-1}} \right)^2 {s_n}^\alpha=\lambda^n 
$$
upon setting $s_n=t_n/\tau$ with $s_0=1$.
By Lemma \ref{recurrencelemma},
we can choose $\lambda=\lambda(\alpha) \in (0,1)$ such that the sequence $\{s_n\}$, and hence $\{t_n\}$,
is increasing and diverges to $+\infty$.
In particular, there is a smallest positive integer $N$ such that $\ell{t_{N+1}}^\alpha \geq 1$.

As in Lemma \ref{decomposition}, define 
$A_n=\{R_n < |z| < R_{n+1}\}$ with
$R_n=1+t_n$, $R_{-1}=1$ ($t_{-1}=0$), and $R_{N+1}=+\infty$ $(n=-1,0,\ldots ,N)$.
Set $\mu_n=\mu \cdot 1_{A_n}$ and $k_n=\Vert \mu_n \Vert_{L_\infty}$. Then
$k_n \leq \ell {t_{n+1}}^\alpha$ by the definition of $\Vert \mu \Vert_{(\alpha)}$.
Applying Lemma \ref{decomposition} yields
\begin{align}\label{S-sum}
(1-|\zeta|)^2|S_{F_\mu|_\mathbb D}(\zeta)| &\leq 12 (1-|\zeta|)^2 \sum_{n=-1}^{N} 
\frac{k_n}{(R_n-|\zeta|)^2}\\ 
&\leq 12 \sum_{n=-1}^{N} \left(\frac{\tau}{\tau+t_n}\right)^2 \cdot \ell {t_{n+1}}^\alpha.
\end{align}
By the recurrence relation \eqref{t-rec}, 
the last term in \eqref{S-sum} coincides with $12\sum_{n=-1}^{N}\lambda^{n+1} \cdot \ell \tau^\alpha$, and moreover
$$
12 \sum_{n=-1}^{N} \lambda^{n+1} \cdot \ell \tau^\alpha <\frac{12 \ell}{1-\lambda} \tau^\alpha
=\frac{12}{1-\lambda} \Vert \mu \Vert_{(\alpha)} (1-|\zeta|)^\alpha.
$$
Therefore,
$$
(1-|\zeta|)^{2-\alpha}|S_{F_\mu|_{\mathbb D}}(\zeta)| \leq C \Vert \mu \Vert_{(\alpha)}
$$
with $C=12/(1-\lambda(\alpha))$.
\end{proof}

\section*{Acknowledgements} 

The author would like to thank Jun Hu for explaining Theorem \ref{junhu} and
for discussions on related problems during his attendance at the conference ICFIDCAA 2024 in Sendai.
He also thanks the referee for providing references related to this work.

\bibliographystyle{amsplain}

\end{document}